\documentclass{article}
\usepackage{amssymb}
\begin{document}
\title{Some algebraic differential equations with few transcendental solutions}
\author{P. X. Gallagher}
\maketitle
\bf Abstract \rm The differential equation $f^{(k)}=f^{(j_1)}...f^{(j_d)}$ with $d > 1$ and each $j_i<k$ has no entire transcendental solutions. In a sense, in almost all cases transcendental meromorphic solutions can also be excluded, and with substantially fewer possible exceptional cases  transcendental solutions are either elliptic or of the form $f(z) = g(e^{cz})$ with $g$ rational and $c$ constant.
\section{ Introduction }
\rm Work culminating in Eremenko, Liao and Ng [3] shows that for each polynomial $ P(X,Y) \neq 0$, each solution $f,$ meromorphic in $\mathbb{C},$  of the differential equation $P(f^{(k)},f) = 0$ with $k$ in $\mathbb{N}$ is either entire or in the class $W$ consisting of rational functions, elliptic functions and functions  $f$ of the form   $f(z) = g(e^{cz})$ with $g$ rational and $c$ constant. Among the tools used in [3] are a refinement of  Wiman-Valiron Theory  due to Bergweiler, Rippon and Stallard [1] and an analysis of the formal Laurent series solutions of certain differential equations.  For a related treatment of another class of autonomous algebraic differential equations using Nevanlinna Theory, see [2].
\\

 We consider meromorphic solutions $f$ (i.e. solutions meromorphic in $\mathbb{C}$) of autonomous algebraic differential equations of the special form
\\
$$f^{(k)} = f^{(j_{1})}...f^{(j_{d})} \hbox{ with  } d > 1 \hbox{ and }  0 \leq j_{1}, ..., j_{d}< k.  \eqno(1) $$
\\
\textsc{Theorem A} \it There are no transcendental meromorphic solutions of (1) with $\infty$ as a Nevanlinna exceptional value. In particular, (1) has no transcendental meromorphic solutions with only finitely many poles.
\\
\\
\rm The second statement, which is the part of Theorem A used in the rest of this paper,  is a special case of Theorem 13.1 in [1], which states that for differential equations of the form $A=0,$ where $A$ is a linear combination with polynomial coefficients of terms of the form $f^{(j_{1})}...f^{(j_{e})}$ with $e \geq 1,$  and assuming $A$ has only one term of top degree in $f$ and its derivatives, there are no transcendental meromorphic solutions with at most finitely many poles. Earlier, Wittich [12] had shown that such equations have no entire transcendental solutions. Recently Zhang and Liao [13] have shown that these equations have no transcendental meromorphic solutions with $\infty$ as a Nevanlinna exceptional value. Our short proof of Theorem A uses the Hayman-Miles Theorem [7].
\\
\\
\textsc{Theorem B} \it Each periodic meromorphic solution of (1) with only finitely many poles in a period strip has the form  $f(z) = g(e^{cz})$ with $g$ rational  and $c$ constant.
\\

\rm The proof of  Theorem B uses  Theorems 2.1 and 2.2 of [1], as in [3]. Also used in the proof of Theorem B, and later, is the fact that the multiplicity $m$ of the pole in any formal Laurent series solution of (1) must satisfy
$$ k = m(d-1) + h, \hbox{ with } h = j_{1} + ...+ j_{d}.  \eqno (2)$$

Key steps in [2] and [3] show that certain autonomous differential equations have at most finitely many formal Laurent series solutions about $z=0.$ This produces a nonzero period for each meromorphic solution with sufficiently many distinct poles. 

We use the special form of (1) to give not finitely many formal Laurent series solutions, but, except in relatively few cases, a special shape for these series solutions still sufficient to  produce a nonzero period for meromorphic solutions with at least two distinct poles. To each differential equation (1) and positive integer $m$ satisfying (2) we associate a polynomial $p.$ The key Lemma 2 relates the roots of $p$ in $\mathbb{N}$ to the shape of formal Laurent series solutions of (1) with a pole of multiplicity $m.$ Using Theorems A and B, this leads to:
\\ 
\\
\textsc{Theorem C} \it Let $f$ be a meromorphic solution of (1) having at least one pole, of multiplicity $m.$ If the associated polynomial $p$ has no positive integer roots, then $f$ is rational. If $p$ has at least one positive integer root, and the greatest common divisor $q$ of these roots satisfies $q > 1,$ then $f$ is in $W.$
\\

\rm The first statement in Theorem C applies to the following cases:
\\
\\
\textsc{Theorem D} \it  If in (1) the $j_i$ are all even or all odd, then, except for the cases $f''= f^{2},$ $f''= f^{3}$, and $f''' = f'^{2},$ all meromorphic solutions of (1) are rational. In the first two exceptional cases there are also elliptic function solutions. In the third case there is a meromorphic solution not in $W:$  an integral of a renormalized Weierstrass $ \wp$-function with period ratio $e^{2 \pi i /6}.$
\\

\rm More detailed versions of Theorems C and D are in sections 3 and 5. 
\\

\rm For the cases of (1) with all $j_{i}=0$ or $1,$ we find all roots of $p$ in $\mathbb{N}.$ Combined with Theorems C and D, and with essential help from the result of Eremenko, Liao and Ng on one case for each odd $k>1$ in which $p$ has two roots in $\mathbb{N}$ and $q=1,$ this leads to:
\\
\\
\textsc{Theorem E} \it Except for $f''' = f'^{2},$ all meromorphic solutions of
$$ f^{(k)} = f^{a}f'^{b} \hbox{ with } a \geq 0, b \geq 0, a+ b > 1 \hbox{ and } k > 1,$$
are in $W.$
\\

\rm Theorem C is largely silent about  cases of (1) for which $p$ has more than one root in $\mathbb{N}.$ However, these cases are relatively rare:

Let $a_{j}$ be the number of $i$ with $j_{i} = j.$ Then (1)  and (2) are equivalent to
$$f^{(k)} = \Pi_{j=0}^{l}(f^{(j)})^{a_{j}},  \hbox{  with  }   k > l \geq 0,   \eqno(3)$$
and
$$ \Sigma_{j=0}^{l}(j+m)a_{j} = k+m. \eqno(4)$$

For positive integers $k,l,m$ with $k > l,$ denote by $A_{k,l,m}$ the set of all ($l+1$)-tuples of nonnegative integers $\bf{a}$ = $(a_{0}, ..., a_{l})$
which satisfy (4).  It turns out that for $\bf{a}$ in $A_{k,l,m}$ the roots $r$ of the corresponding polynomial $p$ in $\mathbb{N}$ satisfy
$m \leq r \leq k+l+2m.$ For fixed $l$ and $m$ each $r \geq k+m$ occurs for only boundedly many $\bf{a}$ in $A_{k,l,m},$ while for (in a sense) almost all $\bf{a}$ in $A_{k,l,m}$ there are no $r < k+m,$ and with even fewer exceptional cases at most one such $r.$ Combined with Theorem C, this leads to:
\\
\\ 
\textsc{Theorem F} \it For positive integers $ k,l,m$ with $l$ and $m$ fixed and $k \to \infty,$
\\
(i) the number of $\bf{a}$ in $A_{k,l,m}$ is  asymptotic to $c_{l,m} k^{l}$ with constants $c_{l,m} > 0;$
\\
(ii)  the number of $\bf{a}$ in $A_{k,l,m}$ for which (3) has a transcendental meromorphic solution having poles of multiplicity $m$ is $O_{l,m}(k^{l-1});$ 
\\
(iii) for $l\geq 2$ and  $\varepsilon > 0,$ the number of $\bf{a}$ in $A_{k,l,m}$ for which (3) has a meromorphic solution not in $W$ having poles of multiplicity $m$ is $O_{l,m,\varepsilon}(k^{l-2+ \varepsilon}).$
\section{ Solutions with few poles}
\textsc{Theorem A} \it  There are no transcendental meromorphic solutions of (1) with $\infty$ as a Nevanlinna exceptional value. In particular, (1) has no transcendental meromorphic solutions with at most finitely many poles. 
\\ 

\rm In the proof of Theorem A we will use some facts from Nevanlinna Theory: To each nonconstant $f$ meromorphic in $\mathbb{C}$ are associated  nonnegative functions of a variable $ r \geq 1:$
 $$ m(f) = m(r,f), N(f) = N(r,f) \hbox{ and } T(f) = m(f) + N(f).$$ 
For definitions of these functions and proofs of the properties (i), (ii) and (iii) below, see [4],[6],[8],[10].
\\
\\
(i) $T(f) \to \infty$  for $ r \to \infty.$
\\
\\
(ii) $ m(f^{d}) = dm(f)$  for $ d $ in $ \mathbb{N},$  and $m(fg) \leq m(f) + m(g).$
\\

For transcendental $f,$
\\
\\
(iii) $m(f'/f) = o(T(f))$  for $ r \to \infty$ off some set  of finite measure.
\\
\\
Combining (iii) with the inequality in (ii) with $g=f'/f$ gives
 $$m(f')  = O(T(f)) \hbox{ off some set of finite measure. }$$
 $N(r,f)$ is a weighted sum, with nonnegative weights, of the multiplicities of the distinct poles of $f$ in $|z| \leq r,$ so
 $$ N(f') \leq 2N(f) = O(T(f)).$$
 It follows that 
\\
\\
(iv) $T(f') = O(T(f))$  off some set of finite measure.
\\
\\
For integers $j$ and $k$ with $0 \leq j <k,$
$$f^{(k)}/f^{(j)} = \Pi_{i=j}^{k-1}(f^{(i+1)}/f^{(i)}), $$
so
$$ m(f^{(k)}/f^{(j)}) \leq  \Sigma_{i=j}^{k-1}m(f^{(i+1)}/f^{(i)}).$$
Combined with (iii) and (iv), this gives
\\
\\
(v) $m(f^{(k)}/f^{(j)}) = o(T(f))$ for $ r \to \infty$ off some set of finite measure.
\\ 

If $f$ has $\infty$ as a \it Nevanlinna exceptional value, \rm meaning
$$N(f) = o(T(f)) \hbox{ for } r \to \infty, \hbox{  i.e. } m(f) \sim T(f) \hbox{ for } r \to \infty,$$ then for each $k$ in $\mathbb{N},$ 
$$ N(f^{(k)}) \leq (k+1)N(f) = o(T(f)) \hbox{ for } r \to \infty , $$
so
\\
\\
(vi) $ T(f^{(k)}) = m(f^{(k)}) + o(T(f))$  for $ r \to \infty.$
\\

In 1989, Hayman and Miles [7] proved the following corrected version of a then 60 year old conjecture of R. Nevanlinna: For each transcendental meromorphic $f,$ and each $k$ in $\mathbb{N},$
 $$T(f) = O(T(f^{(k)}))$$ 
 off some sets of arbitrarily small upper logarithmic density, i.e., for each such $f$ and $k$ and each $\delta > 0,$ there is a constant $B$ and a measurable subset $E$ of $[1, \infty)$ so that $T(f) \leq BT(f^{(k)})$ for all $r>1$ not in $E,$ and 
$$  \int_{E\cap[1,r]} dt/t  < \delta \log r \hbox{ for sufficiently large } r.$$

We will abbreviate \it  sets of arbitrarily small upper logarithmic density \rm by \it arbitrarily small sets. \rm Thus in (v) \it off some set of finite measure \rm can be replaced by \it off some arbitrarily small sets. 

\rm The Hayman-Miles Theorem, combined with (vi), shows that that for transcendental $f$ with $\infty$ as a Nevanlinna exceptional value, and $k$ in $\mathbb{N},$
\\
\\
(vii) $ T(f) = O(m(f^{(k)})),  \hbox{ off some arbitrarily small sets.}$  
\\
\\
\it Proof of Theorem A. \rm For each transcendental meromorphic solution $f$ of (1), 
$$ (f^{(k)})^{d-1} = \Pi_{i=1}^{d}(f^{(k)}/f^{(j_{i})}).$$
Applying (ii), this gives
$$ (d-1)m(f^{(k)}) \leq \Sigma_{i=1}^{d}m(f^{(k)}/f^{(j_{i})}).$$
If $\infty$ is a Nevanlinna exceptional value for $f, $ then this last inequality can be combined with (vii) and (v), giving the contradiction
$$ T(f) = o(T(f)) \hbox{ for  } r \to \infty \hbox{  off some arbitrarily small sets.}$$
\textsc{Theorem B} \it For $d > 1,$ each periodic meromorphic solution of (1) with at most finitely many poles in a period strip has the form $f(z) = g(e^{cz})$ with $g$ rational and $c$ constant.
\\

\rm In the proof of Theorem B we will use the refinement of Wiman-Valiron Theory in [1], which applies in particular to functions $g(w)$ analytic in  $|w| \geq r_{0}$ with  an essential singularity at $\infty:$  For such functions $g,$ and sufficiently large $B>0$  let $G$ be a component of the set of  $w$ with  $ |w| > r_{0}$  and $|g(w)| >B .$   Put $ r_{1} =$ inf$ |w|$ for $w \in G.$ For $r > r_{1},$ denote by $w_{r}$ a point on the part of $|w| = r$ in $G$ at which $|g(w)|$ is maximal and put  $M(r) = |g(w_{r})|.$  $M$ is an increasing function, differentiable off a countable set. Put $a(r) = rM'(r)/M(r).$ Then   
$$a(r) \to \infty \hbox{ for } r \to \infty.  \eqno(5)$$  
by Theorem 2.1 of [1]. By Theorem 2.2 and equation (2.10) in [1] there is a subset $E$ of $[r_{1}, \infty)$ of finite logarithmic measure so that for
 each $j \geq 0,$ 
$$ w_{r}^{j}g^{(j)}(w_{r}) \sim a(r)^{j}g(w_{r}), \hbox{ for } r \to \infty \hbox{ off } E. \eqno(6)$$  
\it Proof of Theorem B. \rm If $f,$ meromorphic in $\mathbb{C},$ has period $\omega \neq 0$ and has only finitely many poles in a period strip, then 
$$f(z) = g(w) \hbox{ with } w = e^{cz} \hbox{ and  }  c\omega = 2\pi i$$ 
and $g$ meromorphic in $\mathbb{C} -  0$ with only finitely many poles. It suffices to show that if $f,$ in addition, satisfies (1), then $g$ does not have an essential singularity at $\infty$ or at $0$.

For $f$ and $g$ related as above, induction gives 
$$f^{(j)}(z) = c^{j} \Sigma_{i=1}^{j}c_{ij}w^{i}g^{(i)}(w) \eqno(7)$$
for $j > 0,$ with constants $c_{ij}$ satisfying $c_{jj} = 1.$  For suitable $r_{0},$  $g$ is analytic in $|w| \geq r_{0}.$ If $g$ has an essential singularity at $\infty,$ then for $G,$ $B,$ $r_{1},$ $r > r_{1}$ and $w_{r}$ as above, choose $z_{r}$ so that $w_{r} = e^{cz_{r}}.$  Combining  (5), (6), and (7) gives 
 $$ f^{(j)}(z_{r}) \sim (ca(r))^{j}g(w_{r}), \hbox{ for } r \to \infty \hbox{ off } E, \eqno(8)$$
for $j>0.$ This holds also for $j=0,$ with equality and for all $r > r_{1}.$ 

If, in addition, $f$ satisfies (1), then (8) for $0 \leq j \leq k$ implies 
$$(|c|a(r))^{k-h} \sim  (M(r))^{d-1} \hbox{ for }  r \to \infty \hbox{ off } E, \eqno(9)$$
with $h = j_{1}  + ... + j_{d}.$

By Theorem A, we may suppose $f$ has at least one pole. Let $m$ be the multiplicity of this pole. Using (2), (9)  simplifies to
 $$|c|M'(r)/M(r)^{1 + {1 \over m}} \sim 1/r, \hbox{ for } r \to \infty \hbox{ off } E. $$
Integrating this over the complement of $E$ in $[r_{1}, \infty)$  gives a contradiction: the the integral of $1/r$ is infinite and the integral of the left side is finite. Thus $g$ does not have an essential singularity at $\infty.$ 

A similar argument shows that $g$ does not have an essential singularity at $0.$ 
\section{ Laurent series and meromorphic solutions}
\rm First, solutions in formal Laurent series with a pole of multiplicity $m:$ Let
$$f(z)= \Sigma _{n=0}^{ \infty }c(n)z^{n-m}, \hbox { with }  m \hbox{ in }\mathbb{ N}  \hbox{ and } c(0) \neq 0.\eqno (10)$$
For each nonnegative integer $j,$
$$f^{(j)}(z)=  \Sigma _{n=0}^{ \infty }(n-m)_{j}c(n)z^{n-m-j} = (-m)_{j}c(0)z^{-m-j} + ... \ \ ,$$
with 
$$(x)_{j} = x(x-1)...(x-j+1) \hbox{  for  } j > 0 \hbox{ and } (x)_{j}= 1 \hbox{ for } j=0.$$
 If (10) satisfies (1), then equality of terms of  least degree is equivalent to
$$k+m= \Sigma_{i=1}^{d}(j_{i}+m) \eqno (11)$$
(a form of (2)), and
$$c(0)^{d-1}=(-m)_{k}/(-m)_{j_{1}}...(-m)_{j_{d}}. \eqno (12)$$
Assuming (11) and (12) are satisfied, then (10) satisfies (1) if and only if
$$(n-m)_kc(n)=  \Sigma (n_{1}-m)_{j_{1}}...(n_{d}-m)_{j_{d}}c(n_{1})...c(n_{d})  \hbox{ for } n \hbox{ in } \mathbb{N},  \eqno (13) $$
 the sum taken over all $d$-tuples of nonnegative integers $n_{1}, ..., n_{d}$ with 
$$ n_{1} + ... + n_{d} = n.$$
From (12), the terms involving $c(n)$ on the right side of (13) are 
$$(-m)_{j_{1}}...(n-m)_{j_{i}}...(-m)_{j_{d}}c(0)^{d-1}c(n)$$ 
$$=((-m)_{k}/(-m)_{j_{i}})(n-m)_{j_{i}}c(n) $$
$$=(-1)^{k-{j_{i}}}(k+m-1)_{k-{j_{i}}}(n-m)_{j_{i}}c(n)$$for $i= 1, ..., d.$  Therefore (13) can be written as
$$p(n)c(n)=s(n)  \hbox{ for } n \hbox{ in } \mathbb{N},\eqno (14) $$
where $p$ is the polynomial of degree $k$ given by
$$p(x):=(x-m)_{k} - \Sigma _{i=1}^{d}(-1)^{k-{j_{i}}}(k+m-1)_{k-{j_{i}}}(x-m)_{j_{i}},  \eqno(15)  $$
and $s(n)$ is the sum of those terms on the right side of (13) for which  
$$ n_{1}+...+n_{d} = n \hbox{ with } 0  \leq n_{i} < n \hbox{ for } i=1, ..., d.$$
This proves:
\\
\\
\textsc{Lemma 1} \it  The differential equation (1) is satisfied by (10) if and only if $m$ satisfies (11), $c(0)$ satisfies (12), and $c(n)$ satisfies (14) for each positive integer $n.$
\\

\rm Let $R$ be the set of roots of $p$ in $\mathbb{N}.$  
\\
\\
\textsc{Remark} If $R$ is nonempty there are infinitely many formal Laurent series solutions (10) of (1):
Let $r$ be the largest element of $R,$  choose $c(0)$ satisfying (12), put $c(n) = 0$ for all $n$ in $\mathbb{N}$ with $n<r,$ and
choose $c(r)$ in $\mathbb{C}$ arbitrarily. Then $s(n) = 0$ for all $n$ in $\mathbb{N}$ with $n \leq r,$ so (14) is satisfied for all such $n,$ using
 $p(r)= 0$ for the case $n=r.$ For $n>r,$ the pair $s(n), c(n)$ is determined recursively satisfying (14): first $s(n)$ by the $c(n')$ with $n' < n$, then $c(n)$ by $s(n)$ and (14), using $p(n) \neq 0.$ 
\\
\\
\textsc{Lemma 2} \it If $p$ has no positive integer roots, then the solutions (10) of (1) are given by $f(z)= c(0)z^{-m},$ with $m$ and $c(0)$ satisfying (11) and (12). If $p$ has at least one positive integer root, and $q$ is the greatest common divisor of all such roots, then the coefficients in each solution (10) of (1) satisfy $c(n)=0$ for each positive integer $n$ not divisible by $q.$
\\
\\
 Proof.  \rm From (14), if $n$ in $\mathbb{N}$ satisfies both $c(n) \neq 0$ and $p(n) \neq 0,$ then $s(n) \neq 0.$ From the form of $s(n)$ it then follows that $n$ is a sum of at most $d$ positive integers $n_{i},$ each smaller than $n$ with $c(n_{i}) \neq 0.$ By induction, each $n$ in $\mathbb{N}$ with $c(n) \neq 0$ is either in $R$ or a sum of elements of $R.$  In particular, if $R$ is empty then $c(n) = 0$ for all $n$ in $\mathbb{N},$ while if $R$ is nonempty then $c(n) = 0$ for all $n$ in $\mathbb{N}$ not divisible by the gcd of the elements of $R.$
\\

Here is a more detailed version of Theorem C in the introduction:
\\
\\
\textsc{Theorem C} \it Let $f$ be a meromorphic solution of (1) with at least one pole, of multiplicity $m,$ and let $p$ be the associated polynomial. 

If $p$ has no positive integer root, then $f$ is rational. 

If $f$ is transcendental, then the greatest common divisor $q$ of the positive integer roots of $p$ is either $1, 2, 3, 4$ or $6.$ 
\\
If $q=2,$ then $f$ is either elliptic or of the form $f(z) = g(e^{cz})$ with $g$ rational and $c$ constant.
\\
If $q = 4,$  then $f$ is elliptic with a period ratio $i.$
\\
If $q = 3$ or $6,$ then $f$ is elliptic with a period ratio $e^{2\pi i/6}.$ 
\\
\\
Proof. \rm Let $f$ be a meromorphic solution of (1) with at least one pole, of multiplicity $m,$ and let $p$ be the associated polynomial. If $p$ has no root in $\mathbb{N},$  Lemma 2 implies that $f$ is rational. Thus for $f$ transcendental $p$ has at least one  root in $\mathbb{N}.$ Let $q$ be the gcd of all such roots.  By Lemma 2, the coefficients $c(n)$ of the Laurent expansion of $f$ about each pole $z_{0}$ satisfy $c(n) = 0$ for all $n$ not divisible by $q.$ It follows that for small $z-z_{0}$ and each $q$-th root of unity $\zeta,$ 
$$f(\zeta(z-z_{0}) + z_{0})  = \Sigma_{n=0}^{\infty}c(n)(\zeta(z-z_{0}))^{n-m} = \zeta^{-m}f(z).$$
By analytic continuation, $$f(\zeta( z  - z_{0}) + z_{0}) = \zeta^{-m}f(z) \hbox{ for all } z.$$
By Theorem A, $f$ has infinitely many poles, in particular at least two distinct poles $z_{0}$ and $z_{1}.$
Since $z_{1}$ also has multiplicity $m,$ it follows that
 $$f(\zeta( z  - z_{0})+ z_{0}) = f(\zeta( z- z_{1})  + z_{1}) \hbox{  for all } z,$$
showing that for each $q$-th root of unity $\zeta \neq 1,$ $$(1-\zeta)(z_{1} - z_{0})  \hbox{ is a nonzero period of } f.$$

For  $q=2,$ the choice $\zeta = -1$ shows that $f$ has a period $\omega_{1} =2(z_{1}-z_{0}).$ If $f$ has only finitely many poles in the period strip $0 \leq Re((z-z_{0}) / \omega_{1}) < 1,$  then $f$ has the form $f(z) = g(e^{cz})$ with $g$ rational and $c$ constant, by Theorem B. If $f$ has infinitely many poles in this strip, then it has a pole $z_{2}$ with $(z_{2} - z_{0})/\omega_{1}$ nonreal. It follows that $f$ is elliptic with periods $\omega_{1}$ and $\omega_{2} = 2(z_{2} - z_{0}).$

For $q > 2 $  the choices $\zeta = \zeta_{q} = e^{2\pi i/q}$ and $\zeta = \bar{ \zeta_{q}}$ give periods 
$$ (1- \zeta_{q})(z_{0} - z_{1}) \hbox{  and } (1- \bar{\zeta_{q}})(z_{1} - z_{0}), \hbox{ with ratio } \zeta_{q}.$$
Since $1+\zeta_{3} = \zeta_{6},$ period ratio $\zeta_{3}$ implies period ratio $\zeta_{6}.$

Since $\zeta_{q}$ is a root of the $q$-th cyclotomic polynomial $\Phi_{q}(x),$ which is irreducible over $\mathbb{Q}$ and has degree $\varphi(q)$ where $\varphi$ is the Euler function [9, p. 279] or [11], it follows that the $ \zeta_{q}^{j}$ for $ 0 \leq j < \varphi (q)$ are linearly independent over $\mathbb{Q}.$  Therefore, so are the sums  $ 1 + \zeta_{q} + ... +  \zeta_{q}^{j-1}$   for $ 1 \leq j \leq \varphi(q).$
For $q > 1,$  multiplying each of these sums by $1 - \zeta_{q}$ shows that the $1- \zeta_{q}^{j}$ for $1 \leq j \leq \varphi(q)$ are linearly independent over $\mathbb{Q}.$ Thus the corresponding $(1-\zeta_{q}^{j} )(z_{1} - z_{0})$ generate a subgroup of rank $\varphi(q)$ in the additive group of periods of $f,$  so the  group of periods of $f$ has rank $ \geq \varphi(q).$ Since $f$ can have at most two independent periods, it follows that $\varphi(q) \leq 2,$ from which $q = 1, 2, 3, 4$ or $6.$
\section{Bounds on the roots of p in $ \mathbb{N}$  }

\textsc{Lemma 3} \it  Each root $r$ of $p$ in $\mathbb{N}$ satisfies $r \geq m,$ with $r=m$ if and only if no $j_{i}$ is $0.$ If $m$ is the only root of $p$ in $\mathbb{N},$ then each meromorphic solution of (1) with poles of multiplicity $m$ is rational. 
\\
\\
Proof. \rm If $r$ is a root of $p$ in $\mathbb{N}$ with $r < m,$ then (15) gives
$$(k+m-r-1)_{k} = \Sigma_{i=1}^{d}(j_{i}+m-r-1)_{j_{i}}(k+m-1)_{k-j_{i}}.   \eqno(16)$$
The left side divided by the $i$-th term on the right simplifies to
$$(k+m-r-1)...(j_{i} + m -r)/(k+m-1)...(j_{i} + m) < 1,$$
showing that (16) is impossible. 

By (15),
$$p(m) = (0)_{k} - \Sigma_{i=1}^{d}(-1)^{k-j_{i}}(k+m-1)_{k-j_{i}}(0)_{j_{i}}.$$
Since $(0)_{j}$ is $1$ for $j=0$ and $0$ for $j > 0,$ 
$$p(m) = (-1)^{k}(k+m-1)_{k}a_{0},$$ 
where $a_{0}$ is the number of $i$ with $j_{i} = 0,$  so $p(m) = 0$ if and only if $a_{0} = 0.$

Put $g = f'.$ If $a_{0} = 0,$ then $g$ is a meromorphic solution of 
$$ g^{(k-1)} = g^{(j_{1}-1)}...g^{(j_{d}-1)},    \eqno(17)$$
with poles of multiplicity $m+1.$ Let $p_{1}$ be the polynomial associated to $ (17)$ and $m+1.$ By (15), first for $(1)$ and $m$ and then for $(17)$ and $m+1,$
$$p(x) = (x-m)p_{1}(x).$$
Assuming $m$ is the only root of $p$ in $\mathbb{N},$ it follows that $p_{1}$ has no  root in $\mathbb{N}$ greater than $ m,$ so no roots in $\mathbb{N}$ at all. Therefore $g$ is rational, so $f$ is rational.
\\

The second statement in Lemma 3 gives a supplement to Theorem C:
\\
\\
\textsc{Corollary} \it  If the associated polynomial $p$ has at most one root in $\mathbb{N},$ then each meromorphic solution of (1) with poles of multiplicity $m$ is in $W.$
\\
\\
Proof. \rm Let $f$ be a meromorphic solution of (1). If $p$ has no roots in $\mathbb{N},$ then $f$ is rational, by Theorem C. If $p$ has exactly one root $r$ in $\mathbb{N},$ then $q = r.$ If $ r> 1, $ then $f$ is in $W,$ by Theorem C. If $r=1,$ then $f$ is rational, by Lemma 3.
\\
\\
\textsc{Lemma 4} \it Let $k,l,m$ be positive integers, with $k > l,$  and let $j_{1}, ..., j_{d}$ be integers in $[0,l]$ satisfying (2). Then each root $r$ of the associated polynomial $p$ satisfies
$$r \leq k+l+2m,  \eqno(18)$$
with equality if and only if $k-l$ is even and each $j_{i} = l.$  
\\
\\
Proof. \rm Put $s=r-m.$ By Lemma 3, $s \geq 0,$ and by (15) the condition $p(r) = 0$ may be written as
$$(s)_{k} = (k+m-1)_{k-l}\sum_{i=1}^{d}(-1)^{k-j_{i}}(l+m-1)_{l-j_{i}}(s)_{j_{i}},$$
which may be rewritten as
$$(k+m)(s)_{k} = (k+m)_{k-l}\sum_{i=1}^{d}(-1)^{k-j_{i}}(j_{i}+m)(l+m)_{l-j_{i}}(s)_{j_{i}}. \eqno(19)$$

To show that (19) implies (18) and to get the conditions for equality in (18), we may assume that $s \geq k+ l + m.$ For $0 \leq j < l,$ 
$$(l+m)_{l-j}(s)_{j} < (l+m)_{l-j-1}(s)_{j+1}, \eqno(20)$$
since this reduces to $s > 2j + 1 + m,$ which follows from $s \geq 2l + m.$ From (19),(20) and (2),
$$(s)_{k} \leq (k+m)_{k-l}(s)_{l}, \hbox{  i.e. } (s-l)_{k-l} \leq (k+m)_{k-l},$$
so $s=k+l+m$ since $k-l > 0.$ By (20), equality is equivalent to $k-l$ even and each $j_{i} = l.$

\section{The cases with $j_{i}$ all even or all odd}
\rm  First, a diophantine equation with binomial coefficients:
\\
\\
\textsc{Lemma 5} \it Let $ s, s_{1}, ..., s_{d}, r$ be positive integers with $s \leq r$ and each $s_{i} \leq r,$ and with $d >1.$ If both
\\
$$ {r-1 \choose s_{1}-1}+ ...+ {r-1 \choose s_{d}-1} = {r-1 \choose s-1} \eqno (21)$$
and
$$ s_{1}+...+s_{d} \leq s,  \eqno (22) $$ 
then $s_{1}=...=s_{d},$ $s= ds_{1},$ and $ r=(d+1)s_{1}.$  
\\
\\
 \it Proof. \rm We will use the fact that for integers $ k, l, n$ with  $0 \leq k<l \leq n$,
$$ {n \choose k} <,=,> {n \choose l} \hbox {  according as  }  k+l<,=,>n.$$
For each $i,$ (22) and $d > 1$ imply $s_{i} < s, $ and (21) gives $${r-1 \choose s_{i}-1}< {r-1 \choose s-1}.$$ 
Therefore $s_{i} +s -2 < r-1,$ so $s_{i} + s \leq r,$ so
$$ { r \choose s_{i}} \leq {r \choose s}. \eqno (23) $$
Writing (21) as
$$s_{1}{r \choose s_{1}}+...+s_{d}{r \choose s_{d}} = s{r \choose s},$$
and using (22) and (23), it follows that there is equality in (22) and  for each $i$ equality in (23), so $s+s_{i}=r,$ giving $s_{1} = ... = s_{d},$ $ds_{1} = s,$ and $ (d+1)s_{1} = r.$ 
\\
\\
\textsc {Lemma 6} \it If $m$ satisfies (11), and all $j_{i}$ are even and some $j_{i}$ is $0,$ then the polynomial $p$ defined by (15) has no positive integer roots unless $k$ is even and all $j_{i}$ are  $ 0,$   in which case $(d+1)m$ is the only positive integer root of $p.$
\\
\\
\it Proof. \rm Let $r$ be a root of $p$ in $\mathbb{N}.$ By Lemma 3, $r \geq m.$ 

For $m \leq r < k+m$, the equation $p(r)=0$ simplifies to 
$$0= \Sigma_{1}^{d} (-1)^{j_{i}}(k+m-1)_{k-j_{i}}(r-m)_{j_{i}}. $$
This is impossible since the terms are nonnegative by the evenness assumption, and a term with $j_{i} = 0$ is positive.

For $r \geq k+m$,  the evenness assumption implies that $p(r)=0$ can be written
$${r-1 \choose k+m-1}= (-1)^{k}\Sigma_{i=1}^{d} {r-1 \choose j_{i}+m-1}.$$
Thus $k$ is even. By Lemma 5 and (11), all $j_{i}=0$ and $r=(d+1)m.$
\\
\\
\textsc{Theorem D} \it If in (1) the $j_{i}$ are all even, or all odd, then, except for the cases $f'' =  f^{3},$ $f'' = f^{2},$  and $f''' = f'^{2},$ all meromorphic solutions of (1) are rational. The first two cases also have elliptic function solutions: a renormalized Jacobi $dn$-function with period ratio $i$ and a renormalized Weierstrass $\wp$-function with period ratio $e^{2 \pi i/6}.$   The third case has solutions not in $W:$ any integral of an elliptic function solution of the second case.
\\
\\
\it Proof. \rm Consider first the special case in which the $j_{i}$ are all even and some $j_{i}$ is $0.$ Combining Lemma 6 and the first statement in Theorem C shows that if $f$ is a nonrational meromorphic solution of (1), then $k$ is even, all $j_{i}$ are $ 0$ and $r = (d+1)m$ is the only root of $p$ in $\mathbb{N},$ where $m$ is the multiplicity of the poles of $f$. In this case $k = (d-1)m,$ so $r= k+2m,$ from which $r \geq 4.$ By Theorem C, $r = 4$ or $6,$ and $f$ is elliptic. 

After a translation we may suppose that $0$ is a pole. Let $z_{1} \neq 0$ be a pole closest to $0.$ 

If $r=4,$ then $(k,m,d) = (2,1,3).$  In this subcase, the proof of Theorem C shows that $\omega_{1} = (1-i)z_{1}$ and $\omega_{2} = (1+i)z_{1}$ are periods. The square inscribed in the disc $|z| \leq |z_{1}|$ with one vertex at $z_{1}$ serves as period parallelogram, with two simple poles, one at $0$ and one at the vertices, so $f$ is a renormalized Jacobi $dn$-function with period ratio $i$.

If $r=6,$ then $(k,m,d) = (2,2,2)$ or $(4,1,5).$ Here $\omega_{1} = (1- \zeta_{6})z_{1}$ and $\omega_{2} = (1- \bar{\zeta_{6}})z_{1}$ are periods. The rhombus with with vertices $0, \omega_{1}, z_{1}, \omega_{2}$ serves as period parallelogram. It is in the disc $|z| \leq |z_{1}|,$ with poles only at the vertices.

 In the subcase $(2,2,2)$ there is a double pole at the vertices, so $f$ is a renormalized Weierstrass $\wp$-function with period ratio $\zeta_{6}.$

In the subcase $(4,1,5)$ there would be only a simple pole at the vertices, so $f$ would have only one pole per period parallelogram, so this does not occur.

In the general case in which the $j_{i}$ are all even or all odd, let $j$ be the smallest $j_{i}.$ We may suppose $j \geq 1.$ If $f$ is a meromorphic solution of  (1), then $g = f^{(j)}$ satisfies
$$g^{(k-j)} = g^{(j_{1}-j)}...g^{(j_{d}-j)},$$ to which the special case above applies. If $g$ is rational, then $f$ is rational. In the remaining two subcases, the multiplicity of the poles of $g$ is $1$ or $2.$ 

 The first subcase does not occur, since the derivative of a meromorphic function does not have simple poles. In the second subcase $j = 1$ and $g$ is a renormalized $\wp$-function with period ratio $\zeta_{6}.$  Since $g(z) = c(z-z_{0})^{-2} + O(1)$ near each pole $z_{0},$ with no residue term and the same $c$ for each pole, $g$ does have a meromorphic integral $f,$ and $f$ is not in $W: f$ (a renormalized Weierstrass $\zeta$-function) is neither rational nor a rational function of an exponential, because it has too many poles, nor is it elliptic, because it has the same residue $-c$ at each pole.
\section{New notation; Cases with all $j_{i} = 0$ or $1$}
\textsc{Notation} \it  In equation (1), denote by $a_{j}$ the number of $i$ with $j_{i} = j.$ Let $l$ be a positive integer, and assume  $a_{j} = 0$ for $j > l.$ In this notation, (1) and the necessary condition (2) for (1) to have a formal Laurent series solution with a pole of multiplicity $m$ are equivalent to
$$ f^{(k)} = \Pi_{j=0}^{l}(f^{(j)})^{a_{j}} \hbox{ with } k>l,   \eqno(24)$$
and
$$  \Sigma_{j=0}^{l}(j+m)a_{j} = k+m.    \eqno(25)$$  
From (15), the polynomial $p$ associated to $a_{0}, ..., a_{l}$ and $m$ may be written as 
$$p(x) = (x-m)_{k} - (-1)^{k}(k+m-1)_{k-l}a(x),  \eqno(26)$$
with
$$ a(x) =\Sigma_{j=0}^{l} (-1)^{j}a_{j}(l+m-1)_{l-j}(x-m)_{j}.    \eqno(27)$$

\rm For the rest of this section, $l=1.$ First we find the roots of $p$ in $\mathbb{N}$ in this case. 
\\
\\
\textsc{Lemma 7} \it Let $m$ and $k$ be positive integers and let $a$ and $b$ be nonnegative integers satisfying
$$ ma + (m+1)b = k+m,$$
and let $p$ be the polynomial associated to $f^{(k)} = f^{a}f'^{b}$ and $ m.$
Then the only  integral root $r$ of $p$ with  $m \leq r < k+m$ is the integer $r,$ if any, satisfying  
$$b(r+1) = k+m.     \eqno(28)$$
The  only integral root $r$ of $p$ with $m+k \leq r \leq k+2m+1$ is
\\
\\
(a) $ r = k+2m+1,$ for $a = 0,$  $b(m+1) = k+m$ and $ k$ odd, 
\\
(b) $r = k+2m,$ for $am = k+m,$  $ b = 0$ and $k$ even,
\\ 
(c) $ r = k+2m-1, $ for $ am = k-1,$   $b = 1$ and  $k$  odd. 
\\
\\
\it Proof. \rm By (26), the only integral roots $r$ of $p$ with $m \leq r < k+m$ are the integral roots $r$  in that interval of 
$a(x).$ By (27) and (25) with $l=1,$
 $$a(x) =  a m - b(x-m)  = k+m - b(x+1).$$
There is a root of $a(x)$ in $\mathbb{N}$ if and only if $b$ divides $k+m,$ in which case the root $r$ satisfies (28). 

By (26) and (27) with $l=1,$ the integral roots $r$ of $p$ with $k+m \leq r \leq  k+2m+1$ are the integral solutions $r$ in that interval of 
$$(r-m)_{k} = (-1)^{k}(k+m-1)_{k-1}(k+m - b(r+1)).   \eqno(29)$$
For $b = 0,$ (29)  means $k$  even and $(r-m)_{k} = (k+m)_{k},$  i.e. $ r = k+2m.$
\\
For $b = 1,$ (29) means $k$ odd and $(r-m)_{k-1} = ( k+m-1)_{k-1}$, i.e. $r = k+2m-1.$
\\
For $b \geq 2,$ then (29) implies $k$ is odd and
$$ (r-m)_{k} \geq (k+m-1)_{k-1}(k+m+2) > (k+m)_{k},$$ 
from which $r > k+2m,$ so $r = k+2m+1.$  From this, Lemma 4 with $l=1$ gives $b(m+1) = k+m.$
Conversely, with odd $k$ and $b(m+1) = k+m,$ the right side of (29) is 
$$(k+m-1)_{k-1}(k+m)(r-m)/(m+1) = (k+m)_{k}(r-m)/(m+1),$$ 
which, for $r=k+2m+1,$  reduces to $(r-m)_{k}.$
\\
\\
\textsc{Theorem E }\it Except for $f''' = f'^{2},$ all meromorphic solutions of 
$$f^{(k)} =  f^{a}f'^{b} \hbox{ with } k>1, a \geq 0, b \geq 0, a + b > 1  \eqno(30)$$
are in $W.$
\\
\\
\it Proof.  \rm By  Theorem C and the Corollary in Section 3, all meromorphic solutions of (30) are in $W,$ unless $p$ has more than one root in $\mathbb{N}$ and $q=1.$

The cases (a),(b), and (c) of Lemma 7 are disjoint, so $p$ has more than one root $r$  in $\mathbb{N}$ if and only if $p$ has exactly two roots in $\mathbb{N}$, one satisfying (28) and the other satisfying one of of (a), (b) or (c). 

Case (a), with $ a = 0,$ and $k$ odd, is covered by Theorem D, with all meromorphic solutions rational, except for the subcase $f''' = f'^{2}$ with solutions not in $W.$

Case (b), with $b=0$ and $k$ even, is also covered by Theorem D, with all meromorphic solutions rational or elliptic.

In Case (c), with $ b = 1$ and $am = k-1$ with odd $k>1,$ the root given by (28) is $r = k+m-1.$ Denoting the gcd of positive integers $u$ and $v$ by $(u,v),$ we have
$$q = (k+m-1, k+2m-1) = (k+m-1, m) = (k-1,m) = m. $$
 By Theorem C, $f$ is in $W$ except possibly for the subcase $m=1.$  For this last subcase $b=1$, $a = k-1$ with odd $k>1$ and $m=1,$ and equation (30)  integrates to
 $$kf^{(k-1)} =  f^{k} + c,  \hbox{ with }  c \hbox{ constant, }$$
from which $f$ is in $W$ by the theorem of Eremenko, Liao and Ng mentioned at the beginning of the introduction. More generally, this theorem covers  all of Case (b) and, by integration as above, all of Case (c).
\section{Small roots of p in $\mathbb{N}$}
For positive integers $k,l,m,$ with $k>l,$ denote by $A_{k,l,m}$ the set of ($l+1$)-tuples of nonnegative integers $\bf{a}$ = $(a_{0}, ..., a_{l})$  satisfying (25). 
\\

Part (iii) of the following lemma, which is  based on a bound for the divisor function,  implies the  infrequency of $\bf{a}$ in $A_{k,l,m}$ for which $p$ has more than one  integral root $ r < k+m.$ 
\\
\\
\textsc{Lemma 8} \it Let $k,l,m$  be positive integers with $l$ and $m$ fixed. For $k \to \infty,$
\\
(i) the number of $\bf{a}$ in $A_{k,l,m}$ is asymptotic to $ k^{l}/(l)!(l+m)_{l+1};$  
\\
(ii) for each integer $r$ in $ [m, k+m),$ the number of $\bf{a}$ in $A_{k,l,m}$ for which $r$ is a root of the associated polynomial $p$ is $O_{l,m}(k^{l-1});$
\\
(iii) for $l \geq 2$ and  $\varepsilon > 0,$  the number of  $\bf{a}$ in $A_{k,l,m}$ for which the associated polynomial $p$ has more than one integral root $r$ in  $[m, k+m)$ is $O_{l,m, \epsilon}(k^{l-2 +\varepsilon}).$
\\
\\
\it Proof. \rm  (i) We show more generally that for each positive integer $l,$ each $(l+1)$-tuple of positive integers $\bf{c} \rm  = (c_{0}, ...,c_{l})$  with $c_{0}$ and $c_{1}$ relatively prime, and each nonnegative integer $k,$ the number $N_{l, \bf{c}}(k)$ of $(l+1)$-tuples of nonnegative integers $(a_{0}, ..., a_{l})$ for which $$c_{0}a_{0} + ...+c_{l}a_{l} = k$$ satisfies
$$ N_{l,\bf{c}}(k) = {k^{l} \over l!c_{0}...c_{l}} + O_{l,\bf{c}}(k^{l-1}  + 1): \eqno(31) $$

For $l=1,$ the relatively prime condition implies both the existence of integers $b_{0}$ and $b_{1}$ with $c_{0}b_{0} + c_{1}b_{1} = k$ and also that the general integral solution of $c_{0}a_{0} + c_{1}a_{1} = k$  is  given by
$$a_{0} = b_{0} + c_{1}t \hbox{  and  } a_{1} = b_{1} - c_{0}t \hbox{ with integral }  t.$$
 Nonnegativity of $a_{0}$ and $a_{1}$ is equivalent to $t$ in the closed interval $[-b_{0}/c_{1}, b_{1}/c_{0}].$ This interval has length $k/c_{0}c_{1},$ so the number of integers in this interval is $k/c_{0}c_{1} + O(1),$ giving (31) for $l=1.$

For $l \geq 2,$ proceed by induction: Writing $c =c_{l},$ $ a =a_{l}$ and $\bf{d}$ = $ (c_{0}, ..., c_{l-1}),$
$$N_{l, \bf{c}}(k) = \Sigma_{0 \leq  a \leq k/c}N_{l-1, \bf{d}}(k-ac)$$
$$= \Sigma_{0 \leq a \leq k/c}({(k-ac)^{l-1} \over (l-1)!c_{0}...c_{l-1}} + O_{l-1,\bf{c}}((k-ac)^{l-2} + 1))$$
$$ = {c^{l-1} \over (l-1)!c_{0}...c_{l-1}}\Sigma_{0 \leq a \leq k/c}({k \over c} -a)^{l-1} + O_{l-1,\bf{c}}(k^{l-1} + 1).$$
For real $x \geq 0,$
$$\Sigma_{0 \leq a \leq x} (x-a)^{l-l} = {x^{l} \over l} + O_{l}(x^{l-1} + 1),$$ 
by comparison with an integral. For $x= k/c,$ this combined with the above calculation gives (31).

For (ii) and (iii) we may assume $k>l.$ Since $(r-m)_{k} = 0$ for all integers $r$ with $m \leq r < k+m,$ it follows that the integral roots of $p$ in this interval are among the integral roots  $r \geq m$ of the polynomial $a(x)$ defined by (30).

(ii) For integral $r \geq m$ the conditions (28) and $a(r) =0$ are two linear equations in $a_{0}, ... , a_{l}.$ Having chosen nonnegative  integers $a_{2}, ..., a_{l}$ arbitrarily  subject to 
$$\Sigma_{j=2}^{l}(j+m)a_{j} \leq k,$$
which may be done in at most $(k+1)^{l-1} = O_{l}(k^{l-1})$ ways, the two conditions reduce to a pair of linear equations for $a_{0}$ and $a_{1}:$
$$ ma_{0} + (m+1)a_{1} = y_{1}$$
and
$$ (l+m-1)_{l}a_{0} - (l+m-1)_{l-1}(r-m)a_{1} = y_{2},$$
with integers $y_{1} $ and $y_{2}$ determined by $k,l,m,r$ and $a_{2}, ..., a_{l}.$
The ratio of the coefficients of $a_{1}$ and $a_{0}$ is positive in the first equation and nonpositive in the second, so there is exactly one 
solution of the pair of equations in rational $a_{0}$ and $a_{1}$, so at most one solution in nonnegative integers $a_{0}$ and $a_{1}.$   

(iii)  From (27) and (25),  $$ a(-1) = \Sigma_{j=0}^{l} a_{j}(l+m-1)_{l-j}(j+m)_{j} $$
$$=(l+m-1)_{l-1}\Sigma_{j=0}^{l} (j+m)a_{j}$$
$$ = (l+m-1)_{l-1}(k+m)$$

 Put  $b(x) = a(x-1)$.  Thus 
$$ b(x) = \Sigma_{i=0}^{l}b_{i}x^{i} \hbox{  with  }  b_{0} = (l+m-1)_{l-1}(k+m).$$
From (27), there are integers $\beta_{i,j}$ (depending also on $l$ and $m$) so that
$$ b_{i} = \Sigma_{j=i}^{l} \beta_{i,j}a_{j} \hbox{  with  } \beta_{i,i} \neq 0, \hbox{  for } i=0, ..., l. \eqno(32)$$ 

Let $l \geq 2.$ If $r_{1}$ and $ r_{2}$ are distinct positive integral roots of $a(x),$ then  $r_{1}+1$ and $r_{2} + 1$ are distinct positive  integral roots of $b(x).$ Thus $b(x) = c(x)d(x)$ with 
$$c(x) = (x-r_{1}-1)(x-r_{2}-1)  = \Sigma_{i=0}^{2}c_{i}x^{i} \hbox{ and } d(x) = \Sigma_{i=0}^{l-2}d_{i}x^{i}.$$
Here $ c_{0} = (r_{1}+1)(r_{2}+1),$ $c_{1} = -(r_{1} + r_{2} + 2),$ $c_{2} = 1$ and  $d_{0}, ..., d_{l-2}$  are integers satisfying 
$$ b_{i} = \Sigma_{j} c_{i-j}d_{j}, \hbox{ for } i = 0, ..., l. \eqno(33)$$ 

In particular, $$( r_{1}+1)(r_{2}+1)d_{0} = (l+m-1)_{l-1}(k+m).$$
 A classical result of Wigert [5, Theorem 315] states that for each positive $\varepsilon$ the number of positive integer divisors of a  positive integer $n$ is $O_{\varepsilon}(n^{\varepsilon}).$
It follows that for each positive $ \varepsilon$  there are  $O_{l,m, \varepsilon}( k^{\varepsilon})$ choices for the triple $r_{1}+1,  r_{2}+1, d_{0},$ so $O_{l,m, \varepsilon}(k^{\varepsilon})$ choices  for the quadruple $c_{0}, c_{1} , c_{2}, d_{0}.$

The number of choices for an $(l-2)$-tuple of nonnegative integers $a_{j}$ with $2 < j \leq l$ satisfying 
$$\Sigma_{j} (j+m)a_{j} \leq k+m.$$   
is  $O_{l,m}(k^{l-2}).$  By (32), each such choice determines an $(l-2)$-tuple of integers $b_{i}$ with $2 < j \leq l.$   Using (33), these $b_{i}$, together with the $c_{i}$ determine the $ d_{i} $ with $0 < i \leq l-2 .$ Using (33) again,  the $c_{i}$ and the $d_{i},$ including $d_{0},$ determine the remaining $b_{i}$  with $0 \leq i \leq 2.$ By (32), all the  $b_{i}$ determine the remaining $a_{j}$ with $0 \leq j \leq 2.$ 
\section{Large roots of p in $\mathbb{N};$ Proof of Theorem F}
Lemma 10 below implies the  infrequency of  $ \bf a $ in $A_{k,l,m}$ for which $p$ has any integral root $r \geq k+m.$
\\
\\
\textsc{Lemma 9} \it For integers $l \geq 1,$ $  s \geq 2,$ $  c \geq 1$  and any integer $t,$ the number of ($l+1$)-tuples of integers $b_{0}, ..., b_{l}$ satisfying both
$$ \Sigma_{i=0}^{l}b_{i}s^{i} = t$$
and
$$ |b_{i}| \leq cs  \hbox{  for  }  i = 0, ..., l$$
is $O_{l,c}(1).$
\\
\\
Proof. \rm It is enough to show that for integers $l \geq 1,$ $ s \geq 2,$ $ c \geq 1,$ the number of ($l+1$)-tuples of integers $c_{0}, ... ,c_{l}$ satisfying both
$$ \Sigma_{i=0}^{l} c_{i}s^{i}=0  \eqno(34) $$ 
and
$$ |c_{i}| \leq 2cs, \hbox{  for  }  i = i, ..., l \eqno(35)$$
is $O_{l,c}(1).$

For each ($l+1$)-tuple of integers $c_{0}, ..., c_{l}$ satisfying (34) and (35), there are integers $d_{0}, ..., d_{l+1}$ with $d_{0}=0$ so that
$$ c_{i} + d_{i} = d_{i+1}s   \eqno(36)$$
for $i = 0, ..., l,$ and
$$ |d_{i}| \leq 4c      \eqno(37)$$ 
for $i = 1, ..., l+1.$ In fact, (34) shows that $c_{0} = d_{1}s$ for an integer $d_{1},$ giving (36) for $i=0,$ and (35) for $i=0$ gives $|d_{1}| \leq 2c.$
Assume $d_{1}, ..., d_{j}$ have been defined  for some $j$ with $1 \leq j \leq l,$  satisfying (36) for $i = 0, ..., j-1$ and (37) for $i= 1, ..., j.$  Then
$$ \Sigma_{i=0}^{j-1}c_{i}s^{i} = \Sigma_{i=0}^{j-1}d_{i+1}s^{i+1} - \Sigma_{i=0}^{j-1}d_{i}s^{i} = d_{j}s^{j}.$$
Combined with (34), this gives (36) for $i=j$ for an integer $d_{j+1}.$ Combined with (35) and (37) for $i=j,$ this gives 
$$|d_{j+1}| = |c_{j} + d_{j}|/s \leq 2c + 4c/s \leq 4c.$$

The number of ($l+2$)-tuples of integers $d_{0}, ... d_{l+1}$ with $d_{0} = 0$ satisfying (37) is $(8c+1)^{l+1} = O_{l,c}(1).$  It follows from (36) that there are   $O_{l,c}(1)$ choices for integers $c_{0}, ..., c_{l}$ satisfying (34) and (35).  
\\
\\
\textsc{Lemma  10} \it For integers $ k \geq 2,$ $ l\geq 1,$ $ m \geq 1,$ and each integer $r\geq k+ m, $ the number of $\bf{a}$ in $A_{k,l,m}$  for which $p(r) = 0$ is $O_{l,m}(1).$ 
\\
\\
Proof. \rm Let $s=r-m.$ Then $s \geq 2$ since $k \geq 2.$ Using (26) and (27), the condition $p(r) = 0$ may be written  
$$ \Sigma_{j=0}^{l}\gamma_{j}a_{j}(s)_{j} = t,$$
with $\gamma_{j} = (-1)^{j}(l+m-1)_{l-j}$ and $t = (-1)^{k}(s)_k/(k+m-1)_{k-l}.$ The $\gamma_{j}$ are nonzero integers, depending also on $l$ and $m,$ and $t$ is an integer depending on $k,l,m$ and $s.$ Since
$$(s)_{j} = \Sigma_{i=0}^{j}\gamma_{i,j}s^{i},$$
with integers $\gamma_{i,j},$ the condition $p(r) = 0$ may therefore be written
$$ \Sigma_{i=0}^{l}b_{i}s^{i} = t, \hbox{   with   }  b_{i} = \Sigma_{j=i}^{l}\gamma_{i,j}\gamma_{j}a_{j}.$$
From (25), $a_{0} + ... + a_{l} \leq (k+m)/m \leq k+1 \leq 2k \leq 2s,$ so 
$$|b_{i}| \leq cs, \hbox{ with } c = 2\hbox{max}|\gamma_{i,j}\gamma_{j}| \hbox{ over } 0 \leq i \leq j \leq l.$$
The $\gamma_{i,j}$ are independent of the other parameters, so $c = O_{l,m}(1).$ Thus  the hypotheses of Lemma 9  are satisfied, giving  $O_{l,m}(1)$ choices for the ($l+1$)-tuple $b_{0}, ... , b_{l}.$ Since each $\gamma_{j,j} = 1$ and each $\gamma_{j}$ is nonzero, this gives $O_{l.m}(1)$ choices for the ($l+1$)-tuple $a_{0}, ..., a_{l}.$ 
\\
\\
\textsc{Theorem F} \it For positive integers $k,l,m$ with $l$ and $m$ fixed and $k \to \infty,$
\\
(i) the number of $ \bf{a} $ in $A_{k,l,m}$ is asymptotic to $c_{l,m}k^{l}$ with constants $c_{l,m} > 0;$
\\
(ii) the number of $\bf{a}$ in $A_{k,l,m}$ for which (3) has a nonrational meromorphic solution having poles of multiplicity $m$ is $O_{l,m}(k^{l-1});$
\\
(iii) for $l \geq 2$ and $\varepsilon > 0, $ the number of $\mathbf{a}$ in $A_{k,l,m}$ for which (3) has a meromorphic solution not in $W$ having poles of multiplicity $m$ is $O_{l,m, \varepsilon}(k^{l-2+ \varepsilon}).$
\\
\\
\it Proof. \rm (i) This follows from  Lemma 8(i).

(ii and iii) For $k,l,m$ positive integers and $\bf{a}$ in $ A_{k,l,m}$, let $p$ be the associated polynomial. For $k > l,$ which we may assume, all roots of $p$ in $\mathbb{N}$  are in the interval $[m,k+l+2m],$ by Lemmas 3 and 4. For each integer $r$ in $[k+m, k+l+2m],$ the number of $\bf{a}$ in $A_{k,l,m}$ for which $r$ is a root of $p$ is $O_{l,m}(1),$ by Lemma 10. It follows that the number of $\bf{a}$ in $A_{k,l,m}$ for which $p$ has at least one root  in $\mathbb{N}$ not in $[m,k+m)$ is $O_{l,m}(1).$

 For $l \geq 2,$ the number of $\bf{a}$ in $A_{k,l,m}$ for which $p$ has more than one integral root in $[m,k+m)$ is $O_{l,m}(k^{l-2+\varepsilon}),$ by Lemma 8(iii). For $l=1$ there are no such $\bf{a},$ since an integral root of $p$ in this interval is a root of the polynomial $a$ of degree $1.$

Combining the above statements, it follows that for $\bf{a}$ in $A_{k,l,m},$ with at most $O_{m}(1)$ exceptions if $l=1$ and
 $O_{l,m}(k^{l-2+\varepsilon})$ exceptions if $l \geq 2,$ the polynomial $p$ has at most one root $r$ in $\mathbb{N}.$ By Theorem C, if $p$ has no root in $\mathbb{N},$  or if $p$ has a unique root $r$ in $\mathbb{N}$ with $r= 5$ or $r > 6,$ then each meromorphic solution $f$ of (3) with poles of multiplicity $m$ is rational, while if $r = 2,3,4$ or $6$, each such solution is in $W$.  The number of $\bf{a}$ \rm in these last four cases is  $O_{l,m}(k^{l-1}),$ by Lemma 8(ii). If $r=1$ is the only root of $p$ in $\mathbb{N},$ then $f$ is rational, by Lemma 3 with $m=1.$ 

P. X. Gallagher, Department of Mathematics, Columbia University, New York, N.Y., 10027, USA
\\
pxg@math.columbia.edu
\end{document}